\documentstyle[10pt,amscd,amssymb]{amsart}

\newcommand{\Ext}{\operatorname{Ext}}

\newcommand{\Hom}{\operatorname{Hom}}

\newcommand{\G}{{\mathcal G}}

\newcommand{\sgn}{\operatorname{sgn}}

\newcommand{\ind}{\operatorname{Ind}}
\newcommand{\res}{\operatorname{Res}}
\newcommand{\tl}{\tilde{\lambda}}
\numberwithin{equation}{section}

\theoremstyle{plain}
\newtheorem{theorem}{Theorem}[section]

\newtheorem{lemma}[theorem]{Lemma}
\newtheorem{prop}[theorem]{Proposition}
\newtheorem{conj}[theorem]{Conjecture}

\newtheorem{problem}[theorem]{Problem}

\begin{document}
\title[Irreducible Specht modules]
{\bf Irreducible Specht modules are  signed Young modules.}
\author{\sc David J. Hemmer}

\address
{University of Toledo\\ Department of Mathematics \\
2801 W. Bancroft\\Toledo\\ OH~43606, USA}
\thanks{Research of the author was supported in part by an NSA
Young Investigator's Grant } \email{david.hemmer@@utoledo.edu}

\date{August 2005}

\begin{abstract}
Recently Donkin defined {\it signed Young modules} as a simultaneous generalization of Young and twisted Young modules for the symmetric group. We show that in odd characteristic, if a Specht module $S^\lambda$ is irreducible, then $S^\lambda$ is a signed Young module. Thus the set of irreducible Specht modules coincides with the set of irreducible signed Young modules. This provides evidence for our conjecture that the signed Young modules are precisely the class of indecomposable self-dual modules with Specht filtrations. The theorem is false in characteristic two.
\end{abstract}
\maketitle
\section{Introduction}

Fayers recently determined \cite{FayersirredSpechttypeA} which Specht modules $S^\lambda$ for the symmetric group $\Sigma_d$ are irreducible in characteristic $p \geq 3$, confirming a conjecture of James and Mathas. For $\lambda$ $p$-restricted or $p$-regular the answer was known, so the problem was for $\lambda$ neither $p$-restricted nor $p$-regular. 

It is easy to show that when $p\geq 3$, if $\lambda$ is $p$-regular and $S^\lambda$ is irreducible then $S^\lambda$ is isomorphic to the Young module $Y^\lambda$. Similarly, if $\lambda$ is $p$-restricted and $S^\lambda$ is irreducible then $S^\lambda$ is isomorphic to the twisted Young module $Y^{\lambda'} \otimes \sgn$.  However if $\lambda$ is neither $p$-restricted nor $p$-regular then $S^\lambda$ cannot be a Young or twisted Young module. We will prove that these irreducible Specht modules instead are signed Young modules.

In characteristic two, $\lambda=(2,2)$ is the only partition \cite{JamesMathasirreduSpechtchar2} which is neither 2-restricted nor 2-regular and for which $S^\lambda$ is irreducible. Furthermore the sign representation is trivial so signed Young modules are all Young modules, and one can easily check that $S^{(2,2)}$ is not a Young module, so our main theorem does not hold in characteristic two. Thus for the remainder of the paper we will assume $k$ is a field of characteristic $p\geq 3$, and we emphasize that most of the results are false without this assumption.

We write $\lambda \vdash d$ for
$\lambda=(\lambda_1, \ldots , \lambda_s)$ a partition of $d$. We also write $|\lambda|=d$.  The Young diagram of $\lambda$ is:

$$[\lambda] = \{(i,j) \in {\bf N} \times {\bf N}  \mid j \leq \lambda_i\}.$$
A partition $\lambda$ is {\it $p$-regular} if there is no $i$ such
that $\lambda_i = \lambda_{i+1} = \cdots = \lambda_{i+p-1}>0$. It is
{\it $p$-restricted} if its conjugate partition, denoted $\lambda
'$, is $p$-regular. We write $\unrhd$ for the usual dominance order on partitions. 

The complex simple $\Sigma _d$ modules are the Specht modules
$\{S^\lambda \mid \lambda \vdash d\}$.  Simple $k\Sigma_d$ modules
can be indexed by $p$-restricted partitions or by $p$-regular
partitions. Both
$$\{D^\lambda := S^\lambda / \mbox{rad}(S^\lambda) \mid \lambda
\mbox{ is $p$-regular }\} $$ and
$$\{D_\lambda = \mbox{soc} (S^\lambda)  \mid \lambda \mbox{ is $p$-restricted }\}$$
are complete sets of nonisomorphic simple $k\Sigma_d$ modules. The
two indexings are related by
 $D^\lambda \cong D_{\lambda'} \otimes \mbox{sgn},$ where sgn is the one-dimensional
 signature representation. We also recall that:
 \begin{equation}
 \label{spechttensorsignisdualspecht}
 S^\lambda \otimes \sgn \cong S_{\lambda'}
 \end{equation}
 where $S_\mu$ denotes the dual of the Specht module $S^\mu$. 

For $\lambda \vdash d$ let $\Sigma_\lambda$ be the usual Young subgroup. The set of Young modules $\{Y^\lambda \mid \lambda
\vdash d\}$ is exactly the set of indecomposable summands of the permutation modules $\{M^\mu:=\ind^{\Sigma_d}_{\Sigma_\mu}k \mid \mu \vdash d\}$. The modules $Y^\lambda \otimes \sgn$ are called {\it twisted} Young modules. 
    
Suppose $\alpha=(\alpha_1, \alpha_2, \ldots , \alpha_m)$ and $\beta=(\beta_1, \beta_2, \ldots , \beta_n)$ where $\alpha \vdash a$ and $\beta \vdash b$ and $a+b=d$. Define the \emph{signed permutation module}:
$$M(\alpha \mid \beta) := \ind^{\Sigma_d}_{\Sigma_\alpha \times \Sigma_\beta}k \boxtimes \sgn.$$
Indecomposable summands of signed permutation modules are called {\it signed Young modules.} Donkin has recently \cite{DonkinsignedYoungmodulepaper} classified the isomorphism classes of signed Young modules. They can be indexed by pairs of partitions $\{(\lambda, \mu) \mid |\lambda |  + p  |\mu |=d\}$ and are denoted $Y^{(\lambda \mid p\mu)}$. As the notation suggests,  $Y^{(\lambda \mid p\mu)}$ is a direct summand of $M(\lambda \mid p\mu)$ with multiplicity one.  Notice that summands of $M(\alpha \mid \emptyset)$ are ordinary Young modules and summands of $M(\emptyset \mid \beta)$ are twisted Young modules, so signed Young modules are a simultaneous generalization. They are all self dual with Specht filtrations.

Our main result is motivated by the fact that irreducible Specht modules in the $p$-regular (resp. $p$-restricted) case are easily seen to be Young (resp. twisted Young) modules. 

\begin{prop}
\label{classiccaseirreduSpechtareYoung}

\begin{itemize}
\item[]
\item[(i)]$S^\lambda \cong S_\lambda$ if and only if $S^\lambda$ is irreducible.
\item[(ii)]$S^\lambda \cong Y^\lambda$ if and only if $\lambda$ is $p$-regular and $S^\lambda$ is irreducible.

\item[(iii)]$S^\lambda \cong Y^{\lambda'} \otimes \sgn$ if and only if $\lambda$ is $p$-restricted and $S^\lambda$ is irreducible.
\end{itemize}
\end{prop}
\begin{pf}Irreducible $k\Sigma_d$ modules are self-dual so one implication in (i) is trivial. For $p>2$, Specht modules are indecomposable and $\Hom_{\Sigma_d}(S^\lambda,S^\lambda) \cong k$ \cite[13.17,13.18]{Jamesbook}. Any self-dual module with these properties must be irreducible, as otherwise the obvious map onto the socle would give a contradiction to the Hom condition. This proves (i).

To prove (ii) we will use the adjoint Schur functor ${\mathcal G}$ and some basic facts about the representations of the Schur algebra. Since we will not use this theory again, we will not describe it here, instead just citing the results we need. Suppose $S^\lambda \cong Y^\lambda$. Since $Y^\lambda$ is self-dual then $S^\lambda$ is irreducible by (i). Now $S^\lambda \otimes \sgn \cong S_{\lambda'}$ so $S_{\lambda'} \cong Y^\lambda \otimes \sgn$. Thus:

\begin{eqnarray*}V(\lambda')&\cong &\G(S_{\lambda'})\,\,\, \cite[3.4.2]{HNspechtfilt}\\
& \cong& \G(Y^\lambda \otimes \sgn)\\ 
&\cong& T(\lambda')\,\,\, \cite[3.4.2]{HNspechtfilt}
\end{eqnarray*}
where $V(\lambda')$ is a Weyl module and $T(\lambda')$ a tilting module, thus is self-dual. If $V(\lambda')$ is self-dual then it must be irreducible and $\lambda'$ must be $p$-restricted (see e.g. \cite[p.87]{Mathasbook}), so $\lambda$ is $p$-regular.

Conversely suppose $\lambda$ is $p$-regular and $S^\lambda$ is irreducible, so $S^\lambda \cong D^\lambda$. But $Y^\lambda$ has a Specht filtration with submodule $S^\lambda$ and subquotients $S^\mu$ with $\mu \rhd \lambda$. Thus $[Y^\lambda:D^\lambda]=1$. But $Y^\lambda$ is indecomposable and self-dual so $Y^\lambda = D^\lambda=S^\lambda$.

Since $S^\lambda \otimes \sgn \cong S_{\lambda'}$ and $D^\lambda \otimes \sgn \cong D_{\lambda'}$, part (iii) follows by essentially the same argument.

\end{pf}

Lemma \ref{classiccaseirreduSpechtareYoung} tells us that an irreducible $S^\mu$ with $\mu$ neither $p$-restricted nor $p$-regular is never a Young or twisted Young module. We will show they are {\it signed} Young modules. The lemma is false in characteristic two, for example $S^{(5,1^2)}$ is a direct sum of two simple modules.

\section{Rouquier blocks for symmetric groups.}
We will make considerable use of abacus combinatorics and the description of blocks of $k\Sigma_d$ by residue contents. Both are thoroughly described in the book \cite{JamesKerberbook}. Recently there has been considerable interest in studying certain blocks of $k\Sigma_d$ known as Rouquier blocks. These are blocks with a special, very large, $p$-core relative to their $p$-weight. In particular a Rouquier block of weight $w$  has an abacus display in which the number of beads on runner $i$ exceeds the number on runner $i-1$ by at least $w-1$ for $i=1,2, \ldots, p-1$. We will use only a few simple properties of such blocks, so will not go into further detail here. For a more thorough description of Rouquier blocks we refer the reader to \cite{FayersirredSpechttypeA}. The importance, for us, of Fayers' work is his proof that irreducible Specht modules induce to Rouquier blocks in a nice way. In particular:

\begin{prop}\cite[Lemmas 3.1-3.3]{FayersirredSpechttypeA}
\label{Fayersrouqiuerblockinductiontheorem}
Suppose $\lambda
\vdash d$ and $S^\lambda$ is irreducible. Then there exists $r
\geq d$ and $\mu \vdash r$ such that:
\begin{itemize}
\item[(i)]$S^\mu$ is irreducible.

\item[(ii)]$S^\mu$ lies in a Rouquier block.

\item[(iii)]$\res ^{\Sigma_r}_{\Sigma_d}S^\mu$ has a direct summand which is filtered by copies of $S^\lambda$.

\item[(iv)]$\ind^{\Sigma_r}_{\Sigma_d}S^\lambda$ has a direct summand which is filtered by copies of $S^\mu$.

\end{itemize}
\end{prop}

To apply Prop. \ref{Fayersrouqiuerblockinductiontheorem}, we need to know the modules in parts (iii) and (iv) are actually semisimple. This is easy in characteristic $p>3$ where $\Ext^1_{k\Sigma_d}(S^\lambda,S^\lambda)$ is known to be zero, but this is not known for $p=3$, so we take another approach:

\begin{lemma}
Let $\lambda, \mu \vdash d$. Let $\alpha, \tilde{\alpha} \vdash a$ and  $\beta, \tilde{\beta} \vdash b$ with $d=a+b$. Then:
\begin{itemize}
\item[(i)] \cite[6.3]{DNsurveyrelating} $\Ext^1_{k\Sigma_d}(Y^\lambda, Y^\mu)=0$.
\item[(ii)] $\Ext^1_{k\Sigma_d}(M(\alpha \mid \beta),M(\tilde{\alpha}\mid\tilde{\beta}))=0$.
\item[(iii)]Suppose $|\tau|+p|\sigma|=d$. Then $\Ext^1_{k\Sigma_d}(Y(\tau \mid p\sigma),Y(\tau \mid p\sigma))=0$.
\end{itemize}
\label{lemmaspechtmodulesdon'tselfextend}
\end{lemma}
\begin{pf} To prove (ii) we apply Mackey's theorem to obtain a direct sum of terms of the form $$\Ext^1_{\Sigma_a \times \Sigma_b}(Y^\tau \boxtimes (Y^\sigma \otimes \sgn), Y^\epsilon \boxtimes (Y^\omega \otimes \sgn)).$$ These are all zero by the Kunneth formula and part (i). Then (iii) follows immediately from (ii).
\end{pf}

\begin{lemma}
\label{lemmareductiontoRockblock}
To prove every irreducible Specht module is a signed Young module, it is sufficient to prove it for irreducible Specht modules in Rouquier blocks. 
\end{lemma}
\begin{pf}
Let $S^\lambda$ and $S^\mu$ be as in Prop. \ref{Fayersrouqiuerblockinductiontheorem} and suppose $S^\mu$ is known to be a signed Young module. Then Lemma \ref{lemmaspechtmodulesdon'tselfextend}(iii) guarantees that $\ind^{\Sigma_r}_{\Sigma_d}S^\lambda$ is actually a direct sum of $S^\mu$'s. Since the collection of Young subgroups is closed under conjugation and intersection, it is immediate from Mackey's theorem that the collection of signed permutation modules, and hence of signed Young modules, is closed under induction and restriction to and from Young subgroups. But $S^\lambda$ is a direct summand of $\res_{\Sigma_d}^{\Sigma_r}\ind_{\Sigma_d}^{\Sigma_r} S^\lambda$, hence a direct summand of $\res_{\Sigma_d}^{\Sigma_r} S^\mu$, and hence a signed Young module.
\end{pf}
 So we need to know which Specht modules in Rouquier blocks are irreducible. This is Fayers' main result:

\begin{prop}\cite[4.1, 4.2]{FayersirredSpechttypeA}
\label{FayerswhichSpechtsareirredu}
Let $S^\lambda$ be in a Rouquier block $B$. Let $\lambda$ have $p$-quotient $(\lambda(0), \lambda(1), \ldots , \lambda(p-1))$. Then $S^\lambda$ is irreducible if and only if $\lambda(0)$ is $p$-restricted, $\lambda(p-1)$ is $p$-regular, $\lambda(i)= \emptyset$ for $0<i<p-1$, and both $S^{\lambda(0)}$ and $S^{\lambda(p-1)}$ are irreducible.
\end{prop}

Let $\kappa$ denote the $p$-core corresponding to $B$. It is an easy exercise with the abacus combinatorics  that raising nodes on runners 0 and $p-1$ as described in the previous proposition corresponds to adding vertical and horizontal $p$-hooks respectively to $\kappa$. So Prop. \ref{FayerswhichSpechtsareirredu} can be restated as:

\begin{prop}
\label{restatewhatareirreduspechtinrockblocks}
Let $S^\lambda$ be irreducible in a Rouquier block $B$ with $p$-core $\kappa$ and $p$-weight $w$. Then $\lambda= \tilde{\lambda}+p\mu$ where $\tilde{\lambda}$ is $p$-restricted, $\mu$ is $p$-regular and $S^\mu$ is irreducible. Furthermore $(\tilde{\lambda})'=\kappa'+p\tau$ where $\tau$ is $p$-regular and $S^\tau$ is irreducible.
\end{prop}

\begin{pf}The proposition is immediate from the observation above it about vertical and horizontal $p$-hooks. The $p\mu$ corresponds to the horizontal $p$-hooks and the $p\tau$ to the vertical $p$-hooks. To translate between the two descriptions, note that $\mu=\lambda(p-1)$ and $\tau=\lambda(0)'$.

\end{pf}

Henceforth fix the $p$-core $\kappa\vdash k$. Since $B$ is a Rouquier block with $p$-weight $w$ we know each runner on the abacus display of $\kappa$ has at least $w-1$ beads more than the previous one. Thus there is a Rouquier block with $p$-core $\kappa$ and weight $i$ for $0 \leq i \leq w$ (and perhaps even larger weights depending on the minimum difference in the number of beads on adjacent runners). We will abuse notation and let $e_B$ denote the block idempotent for any of these blocks, i.e. the block of $\Sigma_{k+pi}$ with $p$-core $\kappa$. Multiplication by $e_B$ projects a module onto its direct summand lying in the Rouquier block.
Further we fix $S^\lambda$ irreducible in $B$ as in Prop. \ref{restatewhatareirreduspechtinrockblocks}, so:
\begin{eqnarray}
\label{partitionlabels}
\lambda&=&\tl+p\mu \nonumber\\
(\tl)'&=&\kappa'+p\tau \\
\sigma&:=&\kappa+p\mu.\nonumber
\end{eqnarray}
Then $\sigma$ is $p$-regular, $S^\sigma$ is irreducible and $\lambda=(\sigma'+p\tau)'$. So $\sigma$ is $\lambda$ with the vertical $p$-hooks stripped off and $\tl$ is $\lambda$ with the horizontal $p$-hooks stripped off.

If $\tau=\emptyset$ then $\lambda=\sigma$ is $p$-regular so $S^\lambda \cong Y^\lambda$ by Prop. \ref{classiccaseirreduSpechtareYoung}(i). Thus we can assume $\tau=(\tau_1, \tau_2, \ldots, \tau_s)$ with $\tau_s \neq 0$. 

We need information about some Littlewood Richardson coefficients. See \cite{Jamesbook} for a complete description of the Littlewood-Richardson rule. We will only need a special case. Specifically:

\begin{lemma}
\label{LRcoeflemma}
Let $\rho\vdash k+pv$ be in the Rouquier block with $p$-core $\kappa$ and weight $v$. Let $v+c\leq w$,  so the block with $p$-core $\kappa$ and weight $v+c$ is also a Rouquier block. Then:
$$e_B(\ind_{\Sigma_{k+pv} \times \Sigma_{pc}}^{\Sigma_{k+pv+pc}}(S^\rho \boxtimes \sgn))$$ has a filtration by Specht modules $S^\epsilon$, where the $\epsilon$ which occur can all be obtained from $\rho$ by adding $c$ vertical $p$-hooks, i.e. the abacus displays for $\epsilon$ and $\rho$ agree except on runner $0$.
\end{lemma}
\begin{pf}The Specht module $S^{1^{pc}}$ is isomorphic to $\sgn$, so the multiplicity of $S^\epsilon$ is given by Littlewood-Richardson coefficient $c(\epsilon; \rho, 1^{pc})$. These coefficients are easy to calculate. 

The $1^{pc}$ guarantees that for each $\epsilon$ the coefficient $c(\epsilon; \rho, 1^{pc})$ is either zero or one. The $\epsilon$ for which it is one are exactly those obtained by adding $pc$ nodes to distinct rows of $\rho$. Thus we must determine which of these $\epsilon$ have $p$-core $\kappa$. 

The $p$-core of $\epsilon$ is determined by its residue content (see \cite{JamesKerberbook} for a description of residue content). In this case, for $S^\epsilon$ to lie in $B$, there must be $c$ nodes of each residue added to the diagram of $\rho$.  They key observation is that, since the $p$-weight of $\rho$ is strictly less than $w$, the abacus configuration implies that every addable node of $\rho$ has the same residue, call it $a$. Since we can add at most one node in each row, the only way to get an addable node of residue $a+1\,$(mod p) is to first add the $p-1$ nodes directly above it. But we are required to add $c$ nodes of each residue. We conclude that the only way to add $pc$ nodes to distinct rows of $\rho$ such that there are $c$ nodes of each $p$-residue is by adding  $c$ vertical $p$-hooks. We remark that not every such $\epsilon$ occurs since there is the further requirement that no two of the added vertical $p$-hooks can intersect the same row.
\end{pf}
We will have need for the special case of Lemma \ref{LRcoeflemma} where $\rho$ is $\kappa$ with only horizontal $p$-hooks added on. In this situation only one $\epsilon$ can occur:
\begin{lemma}
\label{firstinduction}
Let $\rho=\kappa+p\mu$ be in Rouquier block of weight $v=|\mu|$ and let $v+c \leq w$ so the block with $p$-core $\kappa$ and weight $v+c$ is also a Rouquier block. Then:
$$e_B(\ind_{\Sigma_{k+pv} \times \Sigma_{pc}}^{\Sigma_{k+pv+pc}}(S^\rho \boxtimes \sgn)) \cong S^{(\rho'+(pc))'}.$$
\end{lemma}

\begin{pf}Notice that $(\rho'+(pc))'$ is just $\rho$ with $pc$ nodes added to the first column. The result follows from Lemma \ref{LRcoeflemma} and the observation that the second column of $\kappa$, and hence of $\rho$, has exactly $p-1$ fewer nodes than the first column. Thus to avoid two vertical $p$-hooks having a node in the same row, the hooks must all be added to the first column.
\end{pf}

\section{Ladders}
In order to prove irreducible Specht modules are actually signed Young modules the last tool we need is a result of James which guarantees each row in the decomposition matrix of $\Sigma_d$ contains a one. Let $\lambda \vdash d$ and recall that $[\lambda]$ is the corresponding Young diagram. For $r \geq 1$ define the {\it rth ladder} to be the set :
$$\{(i,j) \mid i+(p-1)j=p-1+r\} \subset {\bf N} \times {\bf N}.$$ This ladder is the set of points in ${\bf N} \times {\bf N}$ which lie on the line joining $(r,1)$ to $(1, 1+\frac{r-1}{p-1})$. Define the $p$-regularization $\lambda^R$ of $\lambda$ to be the partition whose Young diagram is obtained from that of $\lambda$ by sliding all the nodes of $[\lambda]$ as far as possible up their ladders. Notice that $\lambda^R=\lambda$ if and only if $\lambda$ is $p$-regular. James showed:
\begin{prop}
\cite[Thm. A]{JamesDecompNumbersII}Let $\lambda \vdash d$. Then:
\begin{itemize}
\item[(i)]$\lambda^R \vdash d$ is $p$-regular and $\lambda^R \unrhd \lambda$.
\item[(ii)]$D^{\lambda^R}$ occurs in $S^\lambda$ with multiplicity one.
\item[(iii)] If $D^\mu$ occurs in $S^\lambda$ then $\mu \unrhd \lambda^R$.
\end{itemize}
\label{Jamestheoremaboutoneineveryrow}
\end{prop}

We also need a lemma which will give us information about dominance relations between $\lambda$, $\mu$, $\lambda^R$ and $\mu^R$. For $r=1, 2, \ldots$, define $l_r(\lambda)$ to be the number of nodes of $[\lambda]$ which lie in the $r$th ladder. For $\lambda, \mu \vdash d$ define $\mu \succ \lambda$ if the largest $r$ for which $l_r(\mu) \neq l_r(\lambda)$ satisfies $l_r(\mu)>l_r(\lambda)$. Moving nodes up and down a ladder does not change the ladder numbers, so this is not a total order. However two $p$-regular partitions with exactly the same ladder numbers are equal, thus $\succ$ is a total order on the set of $p$-regular partitions. Fayers showed:

\begin{lemma}\cite[1.2]{FayersReducibleSpecht}
\label{ladderdominancelemma}
Suppose $\lambda, \mu \vdash d$ are $p$-regular. If $\mu \succ \lambda$ then $\lambda \not \!\unrhd \mu$.
\end{lemma}

\section{Main Theorem}
In order to prove irreducible Specht modules are signed Young modules, it remains by Lemma \ref{lemmareductiontoRockblock} only to show that for $S^\lambda$ irreducible in a Rouquier block that $S^\lambda$ is a signed Young module. To do so we employ an idea of Chuang, originally used in \cite{Chuangoriginalbrouepaper}, of successively inducing then truncating to the block with core $\kappa$.
Recall our fixed partition notation from (\ref{partitionlabels}). Start with $S^{\sigma}:=S_0$ and define:

$$S_1:=e_B(\ind^{\Sigma_{k+p|\mu|+p\tau_1}}_{\Sigma_{k+p|\mu|} \times \Sigma_{p\tau_1}}(S^{\sigma} \boxtimes \sgn)).$$
Lemma \ref{firstinduction} describes $S_1$. Now we induce and truncate again. Define:
 
 $$S_2:=e_B(\ind^{\Sigma_{k+p|\mu|+p\tau_1+p\tau_2}}_{\Sigma_{k+p|\mu|+p\tau_1} \times \Sigma_{p\tau_2}}(S_1 \boxtimes \sgn)).$$
 Continue in this way, inducing and truncating until finally we get:
 $$S:=e_B(\ind^{\Sigma_{d}}_{\Sigma_{d-p\tau_s} \times \Sigma_{p\tau_s}}
 (S_{s-1} \boxtimes \sgn)).$$

\begin{prop}
\label{filtrationofS}
The module $S$ is a direct summand of $$\ind^{\Sigma_d}_{\Sigma_{k+p|\mu|} \times \Sigma_{\tau}}S^{\sigma} \boxtimes \sgn.$$ $S$ has a Specht filtration with exactly one copy of $S^\lambda$ and the rest of the form $S^\epsilon$ where $\epsilon^R \succ \lambda^R$. Furthermore the filtration can be chosen so $S^\lambda$ is on top, i.e. there is a surjection from $S$ to $S^\lambda$.
\end{prop}

\begin{pf} Since induction is transitive and each truncation just picks out a direct summand, $S$ is clearly a direct summand. The second part of the lemma follows from Lemma \ref{LRcoeflemma}. If at each stage of the induction you add a single vertical hook with $p\tau_i$ nodes to the $ith$-column you obtain the single copy of $S^\lambda$; i.e. the Littlewood-Richardson coefficient:

$$c(\lambda; \sigma,1^{\tau_1}, 1^{\tau_2}, \ldots , 1^{\tau_s})=1.$$The other $\epsilon$ which occur all come from moving vertical $p$-hooks from $\lambda$ to columns further to the left. This can only shift nodes to higher ladders $l_r$, and so $\epsilon^R \succ \lambda^R$. Also note that $\lambda \rhd \epsilon$ for any such $\epsilon$, so the filtration given by the Littlewood-Richardson rule can be chosen with $S^\lambda$ at the top.
\end{pf}

We can now prove our main theorem:

\begin{theorem}
\label{main}
Let $S^\lambda$ be irreducible in a Rouquier block.  Then $S^\lambda$ is a signed Young module. Consequently all irreducible Specht modules are signed Young modules.
\end{theorem}
\begin{pf}Let $\lambda$, $\kappa$, $\sigma$, $\mu$ and $\tau$ be as in Prop. \ref{restatewhatareirreduspechtinrockblocks}. Since $\sigma$ is $p$-regular, $S^{\sigma}$ is a Young module by Prop. \ref{classiccaseirreduSpechtareYoung}(ii). Thus 
$$\ind^{\Sigma_d}_{\Sigma_{k+p|\mu|} \times \Sigma_{\tau}}(S^{\sigma} \boxtimes \sgn)$$ is a direct sum of signed Young modules, and in particular $S$ is a direct sum of signed Young modules. Now $S^\lambda \cong D^{\lambda^R}$ by Prop. \ref{Jamestheoremaboutoneineveryrow}. Also for $S^\epsilon$ in the filtration of $S$
$$\lambda^R \not\!\unrhd \epsilon^R$$ by Lemma \ref{ladderdominancelemma}. Thus none of the $S^\epsilon$ have $D^{\lambda^R}$ as a composition factor. So $S$ is a self-dual module with exactly one copy of $D^{\lambda^R}$ which appears in the head. Thus $D^{\lambda^R}$ appears in the socle as well, and so must be a direct summand of $S$. Thus $S^\lambda$ is a summand of $S$, so $S^\lambda$ is a signed Young module. 
\end{pf}

\section{Final Remarks}
An important motivation behind this research was to provide evidence for the following conjecture from \cite{HemmerSpechtfiltrationapper}:

\begin{conj}Let $N \in\mod k\Sigma_d$ be indecomposable and self-dual. Suppose $N$ has a Specht (and hence also a dual Specht) filtration. Then $N$ is isomorphic to a signed Young module.
\end{conj}

Self-dual (equivalently irreducible) Specht modules are the most obvious class of indecomposable self-dual modules with Specht filtrations! We remark that  being indecomposable with both a Specht and a dual Specht filtration is not sufficient to guarantee being a signed Young module- there are examples in \cite{HemmerSpechtfiltrationapper} of such modules which are not self-dual, and therefore not signed Young modules.  We have no conjecture on how to parameterize the set of all indecomposable $k\Sigma_d$ modules with both Specht and dual Specht filtrations.

We now know that irreducible Specht modules are signed Young modules, and Donkin has parameterized the isomorphism classes of signed Young modules. Unfortunately we cannot answer the obvious question: 

\begin{problem}
\label{problemoflabels}
Suppose $S^\lambda$ is irreducible. For what $\mu$, $\tau$ is $S^\lambda \cong Y(\mu \mid p\tau)$? 
\end{problem}

We can handle a few special cases of Problem \ref{problemoflabels}. If $\lambda$ is $p$-regular and $S^\lambda$ is irreducible then $S^\lambda \cong Y^\lambda = Y(\lambda| \emptyset)$. Twisted Young modules are of course signed Young modules, but the labelling is a little trickier. It can be seen easily shown from recent work of Brundan and Kujawa \cite{BrundanKujawa} on the Schur superalgebra that:

\begin{prop}\cite{Kujawacomm} Let $\lambda=\tau+p\mu$ with $\tau$ $p$-restricted. Let $m(\tau)$ be the Mullineux conjugate, i.e. $D_\tau \otimes \sgn \cong D_{m(\tau)}$. Then:
$$Y^\lambda \otimes \sgn = Y(m(\tau) \mid p\mu).$$
\label{kujawalabel}
\end{prop}
Thus Prop. \ref{kujawalabel} together with Prop. \ref{classiccaseirreduSpechtareYoung}(iii) gives the label as a signed Young module for an irreducible $S^\lambda$ when $\lambda$ is $p$-restricted. The open case then is the same that was until recently open for classifying irreducible Specht modules, namely when $\lambda$ is neither $p$-regular nor $p$-restricted.

\bibliography{references}
\bibliographystyle{plain}
\end{document}